\documentclass{article}
\usepackage{amsmath}
\usepackage{amssymb}
\usepackage{color}

\usepackage{latexsym}
\usepackage[english]{babel}
\usepackage[dvips]{graphicx}
\usepackage{verbatim}
\newcommand{\bbibitem}{\bibitem}

\newcommand{\llabel}[1]{{\label{#1}}}

\newcommand{\ffoot}[1]{}

\renewcommand{\r}[1]{(\ref{#1})}

\newcommand{\ex}[1]{} 

\font\tenmsb=msbm10
\font\sevenmsb=msbm7
\font\fivemsb=msbm5
\newfam\msbfam
\textfont\msbfam=\tenmsb
\scriptfont\msbfam=\sevenmsb
\scriptscriptfont\msbfam=\fivemsb
\def\Bbb#1{{\fam\msbfam\relax#1}}

\newcommand{\bi}{\begin{itemize}}
\newcommand{\ei}{\end{itemize}}
\newcommand{\bd}{\begin{description}}
\newcommand{\ed}{\end{description}}
\renewcommand{\i}{\item}

\newcommand{\bqn}{\begin{eqnarray}}
\newcommand{\eqn}{\end{eqnarray}}
\newcommand{\eqnn}{\nonumber\end{eqnarray}}
\newcommand{\eqnl}[1]{\llabel{#1}\end{eqnarray}}

\newcommand{\nn}{\nonumber}
\newcommand{\ba}[1]{\begin{array}{#1}}
\newcommand{\ea}{\end{array}}

\newcommand{\mr}{\Bbb{R}}
\newcommand{\R}{\Bbb{R}}
\newcommand{\mc}{\Bbb{C}}
\newcommand{\N}{\Bbb{N}}

\newcommand{\fine}{\end{document}}

\def \trait (#1) (#2) (#3){\vrule width #1pt height #2pt depth #3pt}
\def \qed{\hfill
        \trait (0.1) (6) (0)
        \trait (6) (0.1) (0)
        \kern-6pt   
        \trait (6) (6) (-5.9)
        \trait (0.1) (6) (0)
\medskip}

\newtheorem{ml}{\bf Lemma}
\newtheorem{Theorem}{\bf Theorem}
\newtheorem{mo}{\bf \underline{{\sl Observation}}}
\newtheorem{mcc}{\bf Corollary}
\newtheorem{Definition}{\bf Definition}
\newtheorem{mpr}{\bf Proposition}
\newtheorem{mproperty}{\bf Property}

\newcommand{\bt}{\begin{Theorem}}
\newcommand{\et}{\end{Theorem}}
\newcommand{\bl}{\begin{ml}}
\newcommand{\el}{\end{ml}}
\newcommand{\bo}{\noindent\begin{mo}\rm}
\newcommand{\eo}{\end{mo}}
\newcommand{\bp}{\begin{mpr}}
\newcommand{\ep}{\end{mpr}}
\newcommand{\bc}{\begin{mcc}}
\newcommand{\ec}{\end{mcc}}
\newcommand{\bdeff}{\begin{Definition}}
\newcommand{\edeff}{\end{Definition}}
\newcommand{\bproperty}{\begin{mproperty}}
\newcommand{\eproperty}{\end{mproperty}}

\newtheorem{mrem}{\bf \underline{{\sl Remark}}}
\newcommand{\brem}{\begin{mrem}\rm}
\newcommand{\erem}{\end{mrem}}

\newcommand{\ppotR}[3]
{

\begin{figure}\begin{center}
~\includegraphics[width=#3truecm]{./#1.eps}\\
\caption{#2}
\llabel{#1}
\end{center}
\end{figure}
\noindent$\!\!$}

\newcommand{\lam}{\lambda}
\newcommand{\g}{\gamma}
\newcommand{\al}{\alpha}
\newcommand{\eps}{\varepsilon}

\newcommand{\pro}{$({\cal P})$\ }
\newcommand{\con}{{\cal C}}
\newcommand{\E}{{\cal E}}

\newcommand{\K}{{\cal K}}

\newcommand{\cc}{ constant control }

\newcommand{\roa}{\rho_A}
\newcommand{\rob}{\rho_B}
\newcommand{\D}{{\cal D}}

\newcommand{\lu}{\lam_1}
\newcommand{\ld}{\lam_2}
\newcommand{\lt}{\lam_3}
\renewcommand{\lq}{\lam_4}
\newcommand{\rocc}{\mbox{\Large $\rho$}_{CC}}
\newcommand{\rorc}{\mbox{\Large $\rho$}_{RC}}
\newcommand{\rorr}{\mbox{\Large $\rho$}_{RR}}

\newcommand{\rr}{{\bf (RR)\ }}
\renewcommand{\cc}{{\bf (CC)\ }}
\newcommand{\rc}{{\bf (RC)\ }}

\newcommand{\V}{\mbox{{\bf V}}}

\newcommand{\aaa}{\mbox{{\bf A}}}
\newcommand{\Aa}{{\cal A} \!\!\!\!{\cal A}}

\newcommand{\Dd}{\mathcal{D}}

\newcommand{\Cc}{\mathcal{C}}
\newcommand{\lunga}{\longrightarrow}
\newcommand{\impl}{\Longrightarrow}
\newcommand{\ra}{\rho_A}
\newcommand{\rb}{\rho_B}
\newcommand{\la}{\lambda}

\newcommand{\ca}{\mathcal{K}}
\newcommand{\del}{\delta}

\newcommand{\vett}[2]{\Big(\begin{array}{c} #1 \\ #2 \end{array}\Big)}

\newcommand{\GUES}{{\bf GUES}}

\newcommand{\csi}{\Xi}

\begin{document}
\begin{center} \noindent
{\LARGE{\sl{\bf  Common Polynomial Lyapunov Functions
for Linear Switched Systems   }}}
\end{center}

\vskip 1cm
\begin{center}
Paolo Mason, Ugo Boscain,

{\footnotesize SISSA-ISAS,
Via Beirut 2-4, 34014 Trieste, Italy}

Yacine Chitour

{\footnotesize  Universit\'e Paris XI,
D\'epartement de Math\'ematiques,
F-91405 Orsay, France}
\end{center}

\vspace{2cm} \noindent \rm
\begin{quotation}
\noindent  {\bf Abstract}
In this paper, we consider linear switched systems $\dot x(t)=A_{u(t)} 
x(t)$, $x\in\R^n$, $u\in U$, and the problem of 
asymptotic stability for arbitrary switching functions, uniform with 
respect to switching ({\bf UAS} for 
short). 
We first prove that, given a {\bf UAS} system, it is always possible to 
build a common polynomial Lyapunov function. Then our main result 
is that the degree of that common polynomial Lyapunov 
function is not uniformly bounded over all 
the {\bf UAS} systems. This result 
answers a question raised by Dayawansa and Martin. A generalization to
a class of
piecewise-polynomial Lyapunov functions is given.

\end{quotation}


\vskip 4cm
\begin{center}
PREPRINT SISSA 52/2004/M
\end{center}

\newpage

\ffoot{
DA FARE:
\bi
\i (cercare se era in un libro di open problems)
\i qualche ref da bacciotti e libro daniel
\i verificare se si puo' passare alla norma delle matrici anche con il 
fattore moltiplicativo che non e' bounded away from zero
\i verificare le forme normali
\i si puo' cambiare la norma delle matrici
\ei
}

\section{Introduction}
In recent years, the problem of stability and stabilizability of switched 
systems has attracted increasing attentions 
(see for instance 
\cite{agr,angeli-lib,bacciotti,SIAM,DM,survey,liberzon-book,xu}), and 
still many questions remain unsolved.

In this paper, we address the problem of existence of common polynomial 
Lyapunov functions for linear switched systems.

By a  switched system, we mean a family of continuous--time dynamical
systems and a rule that determines at each time which dynamical system
 is
responsible of the time evolution. More precisely, 
let $\{f_u:~u\in U\}$ (where $U$ is a subset of 
$\R^m$, $m\in\N$) be a finite or infinite 
set of sufficiently regular vector fields  on
a manifold $M$, and consider the family of dynamical systems:
\bqn
\dot x=f_u(x),~~x\in M.
\eqnl{sw-nl}
The rule  is given by assigning the so-called switching
function, i.e. a measurable function $u(.):[0,\infty[\to U\subset 
\R^m$. Here, we consider the situation in which
the switching function is not known a priori and represents some 
phenomenon (e.g. a disturbance) that is not possible to control.
Therefore, the dynamics defined in \r{sw-nl} also fits into the framework
of uncertain systems (cf. for instance \cite{bub}). 

In the sequel, we use the notations $u\in U$ to label a fixed
individual system and $u(.)$ to indicate the switching function.

These kind of systems are sometimes called ``n-modal systems'', 
``dynamical polysystems'', ``polysystems'', ``input systems''. 
The term ``switched system'' is often reserved to situations in which the 
set $U$ is finite. For the purpose of this paper, we only 
require $U$ to be measurable. For a discussion of various issues related 
to switched systems, we refer the reader to 
\cite{bacciotti,survey,liberzon-book}. 

A typical problem for switched systems goes as follows. 
Assume that, for every $u\in U$, the dynamical system $\dot x=f_u(x)$ satisfies 
a given property (P). Then one can investigate conditions under which 
property (P) still holds for $\dot x=f_{u(t)}(x)$, where $u(.)$ is an arbitrary 
switching function. 

\bigskip
\noindent
In \cite{agr,SIAM,DM,hesp}, the case of switched linear systems was 
considered:
\bqn
\dot x(t)=A_{u(t)} x(t),~~x\in\mr^n,~~A_u\in\mr^{n\times n},
\eqnl{sw-lup}
where $n$ is a positive integer and 
$u(.):[0,\infty[\to U$ is a (measurable) switching function. 
For these systems, the problem of asymptotic stability of the origin, 
uniformly with respect to switching functions, was investigated.\\

Next, we set 
$\aaa:=\big{\{} A_u\ :\ u\in~U\,\big{\}}$ and, to simplify the notation, we 
still call
\underline{switching
function} the measurable matrix-valued map $A(.):=A_{u(.)}$. 
In this way, the switching 
system \r{sw-lup} reads:
\bqn
\dot x(t)=A(t)x(t), \mbox{ where } x\in\mr^n, \mbox{ and } 
A(.):[0,\infty[\to 
\aaa\mbox{ is a measurable map}.
\eqnl{sw-l}
In the following, we assume that:
\bd
\i[(H0)] the set $\aaa$ is a compact 
subset of the set of $n\times n$ real matrices.
\ed
Moreover, the set of switching functions, denoted by $\Aa$, is the 
set
of measurable functions  $A(.):[0,\infty[\to \aaa$.
With our assumptions, for every switching function $A(.)$ and initial
condition $x_0\in\R^n$, the corresponding (Carath\'eodory)
solution of \r{sw-l}
is defined for every $t\geq 0$. We use $\phi^{A(.)}_t (x_0)$ 
to denote the flow of \r{sw-l} at time $t\geq 0$ 
corresponding to the switching function $A(.)$ and starting from 
$x_0$.

Let us recall usual notions of stability used for the system \r{sw-l}.

\bdeff 
\llabel{d-UAS}
Consider the switching system \r{sw-l}. We say that the origin is:
\bd
\i[(S)]  \underline{stable}, if for every $A(.)\in \Aa$ and $\eps>0$, 
there 
exists 
$\del>0$ such that  $\|\phi^{A(.)}_t (x_0)\|\leq\eps$ for every $t\geq 0$, 
$\|x_0\|\leq\del$.  
\i[(US)] \underline{uniformly stable}, if it is stable with $\delta$ not 
depending on $A(.)$.
\i[(U)] \underline{unstable} if it is not stable (i.e. if there exists   
$A(.)\in\Aa$ s.t. the system $\dot x(t)=A(t)x(t)$ is unstable as a 
linear 
time-varying system.)
\i[(AS)] \underline{asymptotically stable}, if it is stable and 
\underline{attractive} (i.e. there exists $\del'>0$ so that, 
for every $A(.)\in\Aa$ and $x_0\in\R^n$ with $\|x_0\|\leq\delta'$, we have 
$\lim_{t\to \infty}\|\phi^{A(.)}_t (x_0)\|=0$). 
\i[(UAS)] \underline{uniformly asymptotically stable}  
if it is uniformly stable and if,
for every $\eps'>0$ and $\del'>0$, there exists $T>0$ such that 
for every switching function $A(.)\in\Aa$, $t\geq T$
and $\|x_0\|\leq\del'$, we have
$\|\phi^{A(.)}_t(x_0)\|\leq\eps'$.
\i[(GUES)] \underline{globally uniformly exponentially stable}, if there 
exist positive constants $M$, $\lam$ such that: 
$\|\phi^{A(.)}_t (x_0)\|\leq M~e^{-\la t}\|x_0\|$,  for every $x_0\in\R^n$, 
$t>0$, $A(.)\in\Aa$. 
\ed
\edeff
Due to the fact that the dynamics is linear in the state variable, the 
local and global notions of stability are equivalent. More precisely,
it was proved in \cite{angeli} that, for system \r{sw-l} subject to {\bf H0}, 
the three notions  {\bf AS}, {\bf UAS}, 
{\bf GUES} and 
the notion of attractivity are all 
equivalent (see also \cite{DM,liberzon-book}). 
In addition, if the system is 
unstable, then there 
exists a 
switching function $A(.)\in\Aa$ and an initial condition $x_0$ such that 
$\lim_{t\to \infty}\|\phi^{A(.)}_t (x_0)\|\to\infty$.
In the following, we just refer to the notions of stability, instability 
and \GUES.

\brem
\llabel{r-class}
Since for the stability issue,  a system of type \r{sw-l}, subject to 
{\bf H0}, is uniquely determined by a compact set $\aaa$ of $n\times n$ 
real matrices, we identify $\aaa$ with the 
corresponding system for the rest of the paper. For instance, when 
we say that $\aaa$ is \GUES, we mean that the corresponding system 
of type \r{sw-l} is \GUES.

\noindent 
We will often consider the problem of determining whether a system, belonging to 
a certain class $C$ of systems of type \r{sw-l} subject to {\bf 
H0}, is \GUES\ or not. Notice that fixing such a class of systems means to 
fix a set of compact subsets of $\R^{n\times n}$ i.e.
$C$ can be identified with a subset of 
$\{\aaa\subset\R^{n\times n}:~~\aaa\mbox{ compact}\}.$
\erem

For a system \r{sw-l} subject to {\bf H0}, it is well known that the \GUES\ 
property is a consequence of the existence of a common Lyapunov function.

\bdeff
\llabel{d-CLF}
A \underline{common Lyapunov function} (\underline{CLF} for short)
$V:\R^n\lunga \R^+$, for a switched system $(S)$ of the type \r{sw-l}, is a continuous
function such $V$ is positive definite (i.e.
$V(x)>0$, $\forall x\neq 0$, $V(0)=0$) and $V$ is strictly
decreasing along nonconstant trajectories of $(S)$.
\ffoot{non uso $\con^1$ per che' voglio trattare PPF} 
\edeff
Vice-versa, it is proved  in \cite{DM} that, given a \GUES\ system of the 
type
\r{sw-l} subject to {\bf (H0)}, 
it is always possible to build a $\con^\infty$ common Lyapunov 
function.

\bigskip
Anyway, the problem of finding a CLF or proving the
nonexistence of a CLF is in general a difficult task. Sometimes, it 
is even easier to prove directly that a system is \GUES\ or unstable. 
An example is provided below by 
bidimensional switched systems.

\subsection{ Single-Input Bidimensional Switched Systems}
Consider a \underline{single input} bidimensional system
of the type:
\bqn
\dot x(t)=u(t)Ax(t)+(1-u(t))Bx(t),
\eqnl{sw-l2}
where $x\in~\mr^2$, A and B are two $2\times2$ real Hurwitz 
matrices and $u(.)$ is a measurable function defined on $\R^+$
and taking values in $U$ equal either to $[0,1]$ or $\{0,1\}$. In the 
sequel, \underline{we call $\csi$} the class of bidimensional 
systems of the above form. This class is parameterized by couples of 
$2\times2$  real Hurwitz matrices. 

\brem\label{r-convex}
Whether systems of type \r{sw-l2} are \GUES\ or not is independent
on the specific choice $U=[0,1]$ or $U=\{0,1\}$. In fact, this is a 
particular instance of a more general result stating that 
the stability properties of systems \r{sw-l} subject to {\bf H0} only 
depend on the convex hull of the set $\aaa$, see Proposition 
\ref{p-convex} below 
and Appendix \ref{s-b}.
\erem

In \cite{shorten}, the authors provide a necessary and sufficient 
condition on the pair $(A,B)$ to share a \underline{quadratic} CLF, but
Dayawansa and  Martin showed in \cite{DM} that
there exist \GUES\  linear bidimensional
systems not admitting quadratic CLF. They also posed the problem of 
finding the  minimal degree of a
\underline{polynomial} CLF. 
More precisely, the problem posed by Dayawansa and Martin is the following:\\\\
{\bf Problem P}:  Define $\Xi_{GUES}\subset\csi$ as the set of 
\GUES\ systems of the type \r{sw-l2}. Find the minimal integer $m$ such 
that every system of $\Xi_{GUES}$ admits a polynomial CLF of degree less 
or equal 
than $m$.

\brem
\llabel{r-dm}
In the problem posed by Dayawansa and Martin, it is implicitly assumed that a 
\GUES\ system always admits a polynomial Lyapunov function and one of our 
results (our Theorem~\ref{t-existence} below) indeed confirms that fact. 
\erem
As for the \GUES\ issue, it was completely resolved in \cite{SIAM},
where a  necessary and 
sufficient condition for a system of type \r{sw-l2} to be \GUES\ was found
directly, without looking for a CLF
(see Section \ref{s-bid} and Appendix \ref{a-A} for more details). 
This is a typical 
example in which it is easier to study directly the stability rather than 
looking for a CLF.

\subsection{Sets of Functions Sufficient to Check \GUES\ }
The concept of Lyapunov function is useful for practical purposes when one can prove 
that, for a certain class of systems, if a CLF exists, then it is possible 
to find one of a certain type and possibly as simple as possible (e.g. 
polynomial with a bound on the degree, piecewise quadratic etc.).

More precisely, consider a class $C$ of systems of type \r{sw-l} in $\R^n$
subject to {\bf H0}, in the sense of Remark \ref{r-class}.
One would like to find a class of functions
${\cal{S}}_C$, identified by a finite number of parameters, which is
sufficient to check \GUES\ for systems belonging to $C$ i.e., if a system of $C$ 
admits a CLF, then it admits one in
${\cal{S}}_C$.  
Once such a class of functions is identified, then in order to verify 
\GUES\, one could use numerical algorithms to check (by varying the
parameters) whether a CLF exists (in which case the
system is \GUES) or not (meaning that the system is not \GUES).

For instance, a remarkable result for a given class $C$ of
systems in $\R^n$ could be the following:\\\\ 
{\bf Claim}: there exists a positive integer $m$ (depending on $n$) such 
that, 
whenever a system of $C$ admits a CLF,
then it admits one that is polynomial of degree less than or equal to $m$.
In other words, the class of polynomials of degree at most $m$ is sufficient 
to check \GUES\ for the class $C$.\\\\
If this result were true,  one could use numerical 
algorithm to check, among all polynomial of degree $m$ (varying the 
coefficients), if there is one that is a CLF. 
Unfortunately, \underline{this claim
is not true}, even for the simplest non trivial case of class of systems 
in $\R^2$, namely systems of type $\csi$ (cf. Equation~\r{sw-l2}). 

The next definition formalizes the idea of class of functions 
sufficient to check \GUES.
\bdeff\llabel{d-ssf}
We say that a 
subset ${\cal{S}}$ of $\con^0(\R^n,\R)$ is finitely (or $m$-finitely)
parameterized if there exist $\Omega\subset \R^m$ for a positive integer 
$m$ and a bijective map $\Psi:\Omega\subseteq\R^m\to {\cal{S}}\subset 
\con^0(\R^n,\R)$. A subset ${\cal{S}}$ of $\con^0(\R^n,\R)$ is said to be 
sufficient to check \GUES\ for a class $C$ of systems of
type \r{sw-l} in $\R^n$ (\underline{{\bf SSF}} for short), if 
every \GUES\ system of $C$ admits a CLF in ${\cal{S}}$. 
A subset ${\cal{S}}$ of $\con^0(\R^n,\R)$ is said to be an 
\underline{($m$-parameters) finite set of 
functions} sufficient to check \GUES\ for a class $C$ of systems of
type \r{sw-l} in $\R^n$ (\underline{finite-{\bf SSF}} for short),
if ${\cal{S}}$ is $m$-finitely parameterized and is an {\bf SSF} for $C$.

If a subset ${\cal{S}}$ of $\con^0(\R^n,\R)$ is not finitely parameterizable 
but is an {\bf SSF} for a class of systems $C$, we call ${\cal{S}}$ an
$\infty$-{\bf SSF}. 
\edeff

\brem
In \cite{blanchini}, a  concept similar  to those introduced in the 
previous definition  was provided and it was called ``universal class 
of Lyapunov functions''.
\erem
Using the previous definitions, the results and the problem  formulated in  
\cite{DM}  can be 
rephrased in the following way:
\bd
\i[R1] for systems \r{sw-l} subject to {\bf H0}, the set  
$\con^\infty (\R^n,\R)$ is an $\infty$-{\bf SSF}; 
\i[R2] for linear bidimensional systems of the class $\csi$ 
(cf. Formula \r{sw-l2}), 
the set of 
quadratic functions is not a {\bf SSF}; 
\i[P] let ${\cal P}^m$ be the set of polynomial functions of two variables with 
degree at most
$m$. What is the minimal $m$ such that ${\cal P}^m$  is a finite-{\bf SSF} for the 
linear bidimensional systems of the class $\csi$?  
\ed
\brem
Notice that, to check numerically the existence of a CLF using the concept 
of finite-{\bf SSF}, one needs some regularity properties of the functions of the 
family, with respect to the parameters (at least continuity). Anyway, this 
discussion is 
out of the purpose of this paper.
\erem
\subsection{Main Results}

We first prove that the implicit assumption of Dayawansa and Martin 
(i.e. that a linear \GUES\ switched system always 
admits a polynomial CLF, cf. Remark 
\ref{r-dm}) is correct in  $\R^n$ (and in particular for bidimensional 
systems of type \r{sw-l2}). 
\bt
\llabel{t-existence}
If the origin is a \GUES\ equilibrium for the 
switched system \r{sw-l} subject to {\bf H0}, then there exists 
a polynomial CLF.
\et
The above result can be stated equivalently as follows. \\\\
{\bf Theorem 1bis} {\it The set of polynomials from $\R^n$ to $\R$, is an 
$\infty$-{\bf SSF} for linear 
switched systems  \r{sw-l} subject to {\bf H0}}.\\\\
The proof of Theorem \ref{t-existence} is given in Section 
\ref{s-existence} and the starting point is the construction of a 
homogeneous and convex CLF $W$, 
following the corresponding argument of \cite{DM}.
The main idea is then to seek for a (homogeneous) polynomial $\tilde {W}$
whose level sets approximate, in some suitable sense, those of $W$
and, finally to show that $\tilde {W}$ is also a CLF.\\
A related result has also been obtained in \cite{blanchini}, using an 
intermediate approximation with polyhedral functions (this step is
implicit in our proof) and starting from the case of discrete approximating 
systems (so-called Euler approximating systems).
\\ 
\\
The core of the paper consists of showing that problem {\bf P} does not 
have a solution, i.e. the minimum degree of a polynomial CLF
cannot be uniformly bounded over the set of all \GUES\ systems 
of the form \r{sw-l2}. More precisely, we have the following:

\bt
\llabel{t-nonexistence}
Let $\Xi_{GUES}\subset \csi$ be  the set of all \GUES\ systems of the type 
\r{sw-l2}. 
If $(A,B)$ is a pair of $2\times 2$ real matrices giving rise to a
system of 
$\Xi_{GUES}$, let $m(A,B)$ be the minimum value of the degree
of any polynomial CLF associated to that system.
Then $m(A,B)$  cannot be bounded uniformly over $\Xi_{GUES}$. 
\et
The above result can be stated equivalently as follows. \\\\
{\bf Theorem 2bis} {\it  Let ${\cal P}^m$ the set of polynomial 
functions from $\R^2\to\R$ of degree at most
$m$. Then ${\cal P}^m$ is not a finite-{\bf SSF} for $\csi$.}\\\\
The proof, given in Section \ref{s-nonexistence}, is based on ideas 
developed in \cite{SIAM}, where
necessary and sufficient conditions for \GUES\ of systems 
\r{sw-l2} are provided. 
We build a sequence of \GUES\ systems
corresponding to 
a sequence of pairs of matrices $(A_i,B_i)$, $i\geq 1$. The sequence of systems
is chosen in such a way that the limit system is uniformly stable but not 
attractive. In particular, that limit system admits a nontrivial periodic 
trajectory whose support $\Gamma$ is a $\con^1$ but not 
a $\con^2$ submanifold of the plane.
To each \GUES\ system of the sequence, one considers any polynomial CLF
$V_{\bar A_i,\bar B_i}$ whose degree is at most $m$. We prove that a 
subsequence of $(V_{\bar A_i,\bar B_i})$ 
converges to a non zero polynomial function $V$ (of degree at most $m$) 
which admits $\Gamma$
as a level set. Since $\Gamma$ is not analytic, a contradiction is reached.

\brem
\underline{The result given by Theorem \ref{t-nonexistence} 
generalizes to 
dimensions higher than 2} as follows. Let 
$(A_i,B_i)$, $i\geq 1$, be a sequence of $2\times 2$ matrices such that 
{\bf i)} the corresponding  systems of type \r{sw-l2} are \GUES, 
{\bf ii)} the 
limit is uniformly stable but not attractive. As 
explained above, for this sequence of systems it is not possible to 
build a sequence of  polynomial CLF of uniformly bounded degree.
Consider now 
the sequence of systems in $\R^n$, $n\geq 2$, of the form 
$\dot{\bar x}=u\bar A_i\bar x+(1-u)\bar  B_i\bar x$ corresponding to the 
matrices:
\bqn
{\bar A}_i=\left(\mbox{
\begin{tabular}{c|cccc}
$A_i$&&0&&\\\hline
&$-1$&...&...&$\vdots$\\
$0$&0&$-1$  &...&$\vdots$\\
&$\vdots$&$\vdots $&$\ddots$&$\vdots$\\
&0&...&...&$-1$
\end{tabular}
}
\right),~~~~{\bar B}_i=
\left(\mbox{
\begin{tabular}{c|cccc}
$B_i$&&0&&\\\hline
&$-1$&...&...&$\vdots$\\
$0$&0&$-1$  &...&$\vdots$\\
&$\vdots$&$\vdots $&$\ddots$&$\vdots$\\
&0&...&...&$-1$
\end{tabular}
}
\right).
\eqn
Each system of the sequence is \GUES\ but the limit system is not (it is 
just uniformly stable). Now, if $V_{\bar A_i,\bar B_i}$, $i\geq 1$, are 
the corresponding 
polynomial CLFs, then they cannot be polynomials of uniformly bounded degree 
since this is not true for the restriction of 
$V_{\bar A_i,\bar B_i}$ to the first two variables.
\erem

\brem {\bf (extension to piecewise polynomial functions (PPF)).}
\llabel{r-ppf}
Another class of functions commonly used to check \GUES\ is that of piecewise
quadratic functions or more generally piecewise polynomial functions (PPF for 
short). Here, by a PPF, we mean a continuous function 
$V\in\con^0(\R^n,\R)$
together with a finite number $q$ of cones $K_j$, $1\leq j\leq q$, based at 
zero and partitioning $\R^n$ so that $V$ is a polynomial function of 
degree 
$d_j$ on $K_j$, $1\leq j\leq q$.
We refer to  $m(V):=\max\{q,d_1,...,d_q\}$  as the \underline{total 
degree} of $V$.

It is tempting to state a version of problem {\bf P} by replacing 
polynomial functions of degree at most $m$ with PPFs  of total degree at most $m$.

Again,  the PPF version of problem {\bf P} does not have a solution for $n=2$,
i.e. the minimum total degree of a piecewise polynomial CLF
cannot be uniformly bounded over the set of all \GUES\ system 
of the form \r{sw-l2}. 
The argument is a simple extension of the proof of 
Theorem~\ref{t-nonexistence}
and it is briefly mentioned in 
Remark~\ref{r-ppf-nonexistence}.
\erem

The last results of the paper concern the existence and the
characterization of a finite-{\bf SSF} for systems of the type \r{sw-l} subject to {\bf H0}.

Let us define the \underline{convex semicone} generated by a set $D\subset 
\R^n$ as  the set of points $\lambda x$ with $\lambda>0$ and $x\in 
co(D)$, where $co(D)$ denotes the convex hull of the set $D$. 
With this definition, the point $x=0$ does not belong to the 
convex 
semicone generated by a set $D$, if $0\notin co(D)$. 

First of all, we prove the following (see Appendix \ref{s-b} for the
argument). 
\bp\llabel{p-convex} 
For every compact subset $\aaa$ of
$\R^{n\times n}$ (i.e. verifying {\bf H0}), let $S_{\aaa}$ be the
system of the type \r{sw-l} associated to $\aaa$.  Then, for $\aaa$
and $\aaa'$ verifying {\bf H0} and generating the same convex
semicone, $S_{\aaa}$ is \GUES\ (resp. uniformly stable) if and only if 
$S_{\aaa'}$ is  \GUES\ (resp. uniformly stable). 
\ep 
Based on
converse Lyapunov theorems, one can deduce some trivial existence
results for finite-{\bf SSF}s. For instance, consider a class $C$ of systems of 
type \r{sw-l} in $\R^n$ subject to {\bf H0} and
satisfying the following property: for every $\aaa\in C$, the convex hull
of $\aaa$ is generated by at most $k$ matrices $n\times n$, where $k$ is a positive
integer. 
One can build a finite-{\bf SSF} for the class $C$ as follows. 
If $C$ does not admit any \GUES\ system, the positive definite 
function $x\mapsto \|x\|^2$ will detect the instability of every element
of $C$. In this case, the subset of ${\cal C}^0(\R^n,\R)$ reduced to 
$x\mapsto \|x\|^2$ is a finite-{\bf SSF} for $C$. Otherwise, 
$C$ admits at least one \GUES\ system. Assume first that a 
finite-{\bf SSF} ${\cal S}$ can be provided for the \GUES\ systems of $C$
and let $V\in {\cal S}$ be a CLF for a fixed \GUES\ system of $C$. 
By associating $V$ to every unstable system of $C$, it is clear that
${\cal S}$ becomes a finite-{\bf SSF} for the whole class $C$.

Therefore, we may simply assume that $C$ is made of \GUES\ systems.
Thanks to Proposition~\ref{p-convex}, the class $C$ can be parameterized by
$k$-tuples of $n\times n$ matrices, defined up to their norm. 
In this way, a $k(n^2-1))$-parameters finite-{\bf SSF} is provided for 
the class $C$. For instance, the class $\csi$ of two-dimensional 
systems of type \r{sw-l2} admits a 6-parameters finite-{\bf SSF}.

The above construction is not explicit and therefore
is not useful to check \GUES. Similarly to Lyapunov functions,
it is then clear that the real challenge for finite-{\bf SSF}s concerns
their \underline{explicit} characterization.
For classes of systems of type \r{sw-l} in $\R^n$, $n\geq 3$, 
that issue is completely open in general. In dimension two,
using the necessary and sufficient conditions for \GUES\ given in  
\cite{SIAM}, \underline{we provide an explicit
5-parameters finite-{\bf SSF}} \underline{for $\csi$}. This is the content 
of Section \ref{s-finite-SSF}.

Clearly, $\csi$ can be parameterized by the pairs $(A,B)$ of $2\times
2$ real Hurwitz matrices, where both $A$ and $B$ are defined up 
to their norm. The construction of the explicit finite-{\bf SSF} goes as follows.  
As done previously, we may assume that the pair $(A,B)$ gives rise to a 
\GUES\ system. By
taking advantage of the complete characterization of \GUES\ systems of the 
class  $\csi$ given in \cite{SIAM}, one can explicitly associate to every
\GUES\ pair $(A,B)$ a CLF as explained next. We start by defining,
from $(A,B)$, a pair $(\tilde{A},\tilde{B})$ giving rise to a system
of $\csi$ which is uniformly stable but not attractive. Such a system
admits a closed trajectory whose support $\Gamma$ is a simple Jordan
closed curve (cf. Sections 
\ref{s-bid}, \ref{s-nonexistence}). We then construct a 
homogeneous positive definite
function $V$ whose level set $1$ is $\Gamma$. We finally show that $V$
is a CLF for $(A,B)$.  Since the set of $(\tilde{A},\tilde{B})$ built
from the \GUES\ pairs $(A,B)$ can be parameterized by using five
parameters, we end up with a five-parameters finite-{\bf SSF}.

\subsection{Structure of the Paper}
Section~\ref{s-existence} contains the proof of Theorem~\ref{t-existence}.
In Section~\ref{ss-bid}, we recall the main ideas from 
\cite{SIAM} needed for the rest of the paper. For sake of completeness,
we provide in Appendix \ref{a-A} the full statement of the main 
result of 
\cite{SIAM}.
The proof of Theorem~\ref{t-nonexistence} is given in Section~\ref{s-nonexistence}
and the explicit construction of a five-parameters finite-{\bf SSF} for systems of type 
\r{sw-l2} is provided in Section \ref{s-finite-SSF}.
Finally, in Appendix \ref{s-b}, we prove Proposition \ref{p-convex}. 


\section{Existence of Common Polynomial Lyapunov Functions}
\llabel{s-existence}
In this section, we prove Theorem \ref{t-existence}. The starting point
of the argument follows the first part of the proof of an analogous 
result in \cite{DM}.\\\\
We define the function $V:\R^n\to \R^+$ by: 
$$ V(x)=\sup_{A(.)\in\Aa}
\int_0^{+\infty} \|\phi^{A(.)}_t (x)\|^2 \mathrm{d} t\ .
$$
The function $V$ is well defined since there exist positive constants $C,\,\mu$
such that, for all $t\geq 0$ and $x\in\R^n$:
$$\|\phi^{A(.)}_t (x)\|\leq C\,e^{-\mu t} \|x\|.
$$
Note that $V$ is homogeneous of degree $2$ and continuous. 
In addition, we next show that $V$ is strictly convex. That fact will be crucial
later in the argument.
Fix $x,y\in\R^n$ and $x\not =y$. Let $A(.)$ be a switching function.
The function $x\mapsto \|\phi^{A(.)}_t (x)\|^2$ is strictly convex.
Moreover, for every $\la\in ]0,1[$, by
compactness of $\aaa\,$, the expression:
$$\la\,\|\phi^{A(.)}_t (x)\|^2+(1-\la)\,
\|\phi^{A(.)}_t (y)\|^2-\|\phi^{A(.)}_t (\la x+(1-\la) y)\|^2,$$ 
is nonnegative for every $t\geq 0$ and is
bounded from below by a positive 
constant on some interval $[0,\bar{t}]$, uniformly with respect to $A(.)$.
\\
Therefore, dividing the integration interval into the two intervals 
$[0,\bar{t}]\,$ and $[\bar{t},+\infty ]$ and taking a maximizing sequence of
switching functions for $V(\la\,x+(1-\la)\,y)$, we have:
$$
V(\la\,x+(1-\la)\,y)<\la\,V(x)+(1-\la)\,V(y),\quad \forall \la\in 
]0,1[,\quad 
\forall x,y\in\R^n\,.
$$  
\\ 
It is shown in \cite{DM} that $V$ is a CLF. Nevertheless,
we need to consider at least $\Cc^1$ Lyapunov functions, therefore we 
define:
$$
\tilde{V}(x)=\int_{SO(n)} f(R)\,V(Rx)\,\mathrm{d}R\,.
$$ where 
$f:SO(n)\lunga [0,+\infty[$ is a smooth function with support on a 
small 
neighborhood of the identity matrix and $\int_{SO(n)} f(R)\,\mathrm{d}R=1$ .\\
In \cite{DM}, it is also shown that $\tilde{V}$ is a smooth CLF except at 
the origin. Moreover, since $V$ is homogeneous of degree $2$ and
strictly convex, it follows that $\tilde{V}$ also satisfies such properties. 
\\
We consider now the function $W(x)=\sqrt{\tilde{V}(x)}$ , which is a continuous,
positively homogeneous CLF. Therefore, $W^{-1}(1)$ is a compact set.
Using the fact that the set $\{ x:\,W(x)<1\}$ is 
strictly convex, we construct a polynomial CLF $\tilde{W}$ by 
approximating
the level sets of $W$. For this purpose, we need the following preliminary
result which describes a continuity property of the function $\nabla W(y)
\cdot Dx$ with respect to $x,\,y,\,D$.
\bl
Let us set:
$$
M:=\displaystyle{\min_{x\in W^{-1}(1)\ D\in \aaa} \big[-\nabla 
W(x)\cdot Dx\big]}.
$$
Then, for every $\eps\in (0,M)$, there exists $\del\in (0,1)$ 
such that, for every $x,y\in W^{-1}(1)$ with $\nabla W(y)\cdot x>1-\del$ and every 
$D\in \aaa$, one has: 
$$
\nabla W(y)\cdot Dx<-\eps.
$$
\el
{\bf Proof of the Lemma.} 
First of all, notice that $M$ is well defined since it is the infimum of a 
continuous function over a compact set. Moreover, $M>0$ because $W$ is a 
CLF.\\
Since, by homogeneity,
$\nabla W(y)\cdot y= W(y)=1$, we have:
$$
\nabla W(y)\cdot x=1-\nabla W(y)\cdot (y-x),
$$ 
and then the hypothesis is equivalent to $\nabla W(y)\cdot (y-x)<\del\, $.\\
Reasoning by contradiction, assume that there exists a sequence $(x_j,y_j,D_j)$ such that
$\nabla W(y_j)\cdot D_j x_j\geq -\eps$ and $\nabla W(y_j)\cdot (y_j-x_j)$
converges to $0$ as $j$ goes to infinity. By compactness, we can find a subsequence of 
$(x_j,y_j,D_j)$ converging to $(\bar{x},\bar{y},\bar{D})$ and therefore, by
continuity, $\nabla W(\bar{y})\cdot \bar{D} \bar{x}\geq -\eps$ and 
$\nabla W(\bar{y})\cdot (\bar{y}-\bar{x})=0$.\\
Therefore $\bar{y}-\bar{x}$ belongs to the tangent space at $\bar{y}$
of the strictly convex set $W^{-1}([0,1])$. Since $\bar{x}$ also
belongs to the boundary of that set, it must be $\bar{y}=\bar{x}$. 
It implies $\nabla W(\bar{y})\cdot \bar{D} \bar{x}=\nabla 
W(\bar{x})\cdot \bar{D} \bar{x}\leq -M$ and we reach a contradiction. 
\hspace{\stretch{1}}$\Box$
\brem
Taking $-x$ instead of $x$, one obtains that 
for every $x,y\in W^{-1}(1)$ and every $D\in \aaa$, then
$\nabla W(y)\cdot x<-1+\del\impl\nabla W(y)\cdot Dx>\eps$.
\erem
To conclude the proof of the theorem, we take $\del\in (0,1)$ corresponding to some $\eps$ 
as in the lemma above, and for every $y\in W^{-1}(1)$ we consider the open sets 
$B_y=\big\{ x\in \R^n\,:\nabla W(y)\cdot x>1-\del/2 \big\}$.
Since $y\in B_y$, we have that $\,\{B_y\}_{y\in W^{-1}(1)}\,$  is an open covering of the 
compact set $W^{-1}(1)\,$, and therefore we can find $y_1,\ldots,y_N\,$ points of 
$W^{-1}(1)\,$ such that the union of $B_{y_k}$, $k=1,\ldots ,N$, covers $W^{-1}(1)\,$.\\
Let us define: 
$$
\tilde{W}(x):=\sum_{k=1}^N (\nabla W(y_k)\cdot 
x)^{2p}\,.
$$
We claim that, for an integer $p$ large enough, $\tilde{W}$ is a polynomial CLF.
For $D\in \aaa$ and $x\in \R^n$, $x\not =0$, we have: 
\bqn\label{comp0}
\nabla \tilde{W}(x)\cdot Dx=2p \sum_{k=1}^N \big(\nabla W(y_k)\cdot x\big)
^{2p-1}\nabla W(y_k)\cdot Dx,
\eqn
and we want to show that $\nabla \tilde{W}(x)\cdot Dx<0$. By homogeneity, it is enough to do it
for $x\in  W^{-1}(1)$. 
Set: $$
K:=\displaystyle{\max_{x,y\in W^{-1}(1),\,D\in\aaa}\,\nabla W(y)\cdot Dx}.
$$
If, for some index $k$ in $\{1,...,N\}$, one has $|\nabla W(y_{k})\cdot x|\leq 1-\del$.
Then:
$$
|\big(\nabla W(y_k)\cdot x\big)^{2p-1} \nabla W(y_k)\cdot Dx |\leq (1-\del)^{2p-1}K.
$$
Otherwise, if the inequalities $1-\del/2\geq |\nabla W(y_{k})\cdot x|>1-\del$ hold,
then, by the previous lemma and remark, one has that the corresponding term in the 
summation must be negative.\\
Finally, since by the definition of the points $y_k$, there exist at 
least two distinct indices $k_1$ and $k_2$ such that 
$x\in B_{y_{k_1}}$ and $\,-x\in B_{y_{k_2}}$ we have that:
$$
\big(\nabla W(y_{k_i})\cdot x\big)^{2p-1}\nabla W(y_{k_i})\cdot Dx
<-(1-\del/2)^{2p-1}\eps.
$$
Summing up, we deduce that:
$$\nabla \tilde{W}(x)\cdot Dx<2p \Big(-2(1-\del/2)^{2p-1}\eps+(N-2)(1-\del)
^{2p-1}\,K\Big)=-4p(1-\del/2)^{2p-1}\eps\Big(1-\frac{K(N-2)}{2\eps}
\Big(\frac{1-\del}{1-\del/2}\Big)^{2p-1}\Big).
$$
For $p$ large enough, the right-hand side of previous 
expression is negative, uniformly with respect to $D\in \aaa$ and $x\in  W^{-1}(1)$.
The theorem is proved. 
\brem
One can also check that the level set ${\tilde W}^{-1}(1)$ approximates, as $p$ tends to 
$+\infty$, the corresponding level set of the function $\max_{k=1,\ldots,N} |\nabla W(y_k)
\cdot x|$ (which is a polytope) and, therefore, the latter is a CLF as well (cf. 
\cite{blanchini}).
\erem

\section{Necessary and Sufficient Conditions for \GUES\ of 
Bidimensional Systems}
\llabel{ss-bid}
\llabel{s-bid}
Consider the following property:
\bd
\item{\pro~}  
The bi-dimensional switched system given by:
\bqn
\dot x(t)=u(t)Ax(t)+(1-u(t))Bx(t), \mbox{ where } u(.):[0,\infty[\to 
[0,1],
\eqnl{sw-l2b} 
is \GUES\ 
at the
origin.
\ed
In this section, together with  Appendix \ref{a-A}, 
we recall the main ideas from \cite{SIAM}, 
to get a necessary and sufficient condition on $A$
and $B$ under which
\pro holds, or under which we have at least uniform stability. The full 
statement of the Theorem is reported in Appendix \ref{a-statment}. 
\brem
\llabel{r-liboni}
Recall that, by Proposition 1 (proved in Appendix \r{s-b}), the  
necessary and sufficient condition for stability of the system \r{sw-l2b} 
are the same if  we assume $u(.)$ taking values in $\{0,1\}$ or 
in $[0,1]$, or if we multiply $A$ 
and $B$ by two arbitrary positive 
constants.
\erem
Set $M(u):=uA+(1-u)B$, $u\in[0,1]$. In the class of constant functions the
asymptotic stability of the
origin of the system \r{sw-l2b} occurs if and only if the
matrix $M(u)$ has eigenvalues with strictly negative
real part for each $u\in[0,1]$. So 
this is  a necessary condition for \GUES. On
the other hand it is
know that if $[A,B]=0$ then
the system \r{sw-l2b} is \GUES.
So, in what follows, we  always assume the conditions:
\bd
\item[H1:] 
Let $\lam_1,\lam_2$ 
(resp. $\lam_3,\lam_4$) be  the eigenvalues of $A$ (resp. B). Then
$Re(\lu),$ $Re(\ld),$ $Re(\lt),$ $Re(\lam_4)<0$. 

\item[H2:] $[A,B]\neq0$ (that implies that neither $A$ nor $B$ is
proportional to the identity).
\ed
For simplicity we will also assume:
\bd
\item[H3:] $A$ and $B$ are diagonalizable in $\mc$ (notice that if {\bf 
H2} 
and
{\bf H3}  hold then $\lu\neq\ld$,
$\lt\neq\lam_4$).  
\item[H4:] Let $\V_1,\V_2\in\mc P^1$ 
(resp. $\V_3,\V_4\in \mc P^1$) be  the eigenvectors of $A$ (resp. B).
Then $\V_i\neq\V_j$ for $i\in\{1,2\}$, $j\in\{3,4\}$
(notice that, from {\bf H2} and {\bf H3}, the $V_i$'s are uniquely defined,
$\V_1\neq\V_2$ and $\V_3\neq\V_4$, and {\bf H4} can be violated only when 
both $A$ and $B$ have real eigenvalues). 
\ed 
All the other cases in which {\bf H1} and {\bf H2} hold are the following: 
\bi 
\item 
$A$ or $B$ are not diagonalizable. 
This case (in
which \pro can be true or false) can be treated with techniques entirely
similar to the ones of \cite{SIAM}. 
\item 
$A$ or $B$ are  diagonalizable, but  
one eigenvector of $A$ coincides
with one eigenvector of $B$. In this case, using arguments similar to 
those of \cite{SIAM}, it possible to conclude that \pro is true.  
\ei 
We will call respectively {\bf (CC)}
the case where both matrices have non-real eigenvalues, 
{\bf (RR)} the case where both matrices have real eigenvalues
and finally {\bf (RC)} the case where one matrix has real eigenvalues and
and the other non-real eigenvalues. 

Theorem \ref{stabb}, reported in Appendix \ref{a-statment}, gives necessary
and sufficient conditions for  the
stability of the system \r{sw-l2b} in
terms of  three (coordinates invariant) parameters 
given below in Definition \ref{zappa}. 
The first two parameters, $\roa$ and $\rob$,
depend on the eigenvalues of $A$ and $B$ respectively, 
and the third parameter $\K$ depends on $Tr(AB)$,
which is a Killing-type pseudo-scalar product in the
space of $2\times2$ matrices. As explained in \cite{SIAM}, the parameter 
$\K$ contains the
inter-relation between the two systems $\dot x=Ax$ and $\dot x=Bx$, and it
has a precise geometric meaning. It is in 1--1  correspondence with  the 
cross ratio of the   four points
in the projective line $\mc P^1$ that corresponds to the four 
eigenvectors of $A$ and $B$.

\bdeff
Let $A$ and $B$ be two $2\times2$ real matrices and suppose that 
{\bf H1}, {\bf H2}, 
{\bf H3} and {\bf H4}
hold. Moreover choose the labels $(_1)$ and $(_2)$  (resp.  $(_3)$ and $(_4)$) 
so that $|\ld|>|\lu|$ (resp. $|\lq|>|\lt|$) if they are
real 
or $Im(\ld)<0$ (resp. $Im(\lam_4)<0$) if they are complex.
Define:
\bqn
&&\roa:=-i\frac{\lu+\ld}{\lu-\ld};
~~~~~~\rob:=-i\frac{\lt+\lam_4}{\lt-\lam_4};
~~~~~~\K:=2\frac{Tr(AB)-\frac{1}{2}
Tr(A)Tr(B)}{(\lu-\ld)(\lt-\lam_4)}.\nn
\eqn
Moreover, define the following function of $\roa,\rob,\K$:
\bqn
\D:=\K^2+2\roa\rob\K-(1+\roa^2+\rob^2).
\eqn
\label{zappa}
\edeff
Notice that $\roa$ is a positive real number if and only if 
$A$ has non-real eigenvalues and $\roa\in i\mr$, $\roa/i>1$ if and 
only if $A$ has real eigenvalues. The same holds for $B$. Moreover 
$\D\in\mr$. 
\brem
Under hypotheses {\bf H1} to {\bf H4}, using 
a suitable 3-parameter changes of coordinates, 
it is always possible to put the 
matrices $A$ and $B$, up the their norm, (cf. Remark \ref{r-liboni}) in 
the normal forms given in Appendix \ref{a-normalforms}, where
$\roa,\rob,\K$ appear explicitly  (see \cite{SIAM} for more 
details).
\erem
The parameter $\K$ contains
important information about the matrices $A$ and $B$. They are stated in
the following Proposition that can be easily proved  using the normal 
forms given in Appendix \ref{a-normalforms}.
\bp 
Let $A$ and $B$ be as in Definition 
\ref{zappa}. Then: 
{\bf i)} if $A$ and $B$ have both complex eigenvalues,  then
$\K\in\mr$ and $|\K|>1$;
{\bf ii)} if $A$ and $B$ have both real eigenvalues, then
$\K\in\mr\setminus\{\pm 1\}$;
{\bf iii)} $A$ and $B$ have one complex and the other real eigenvalues 
if and only if $\K\in i\mr$.
\label{vanga}
\ep
\ppotR{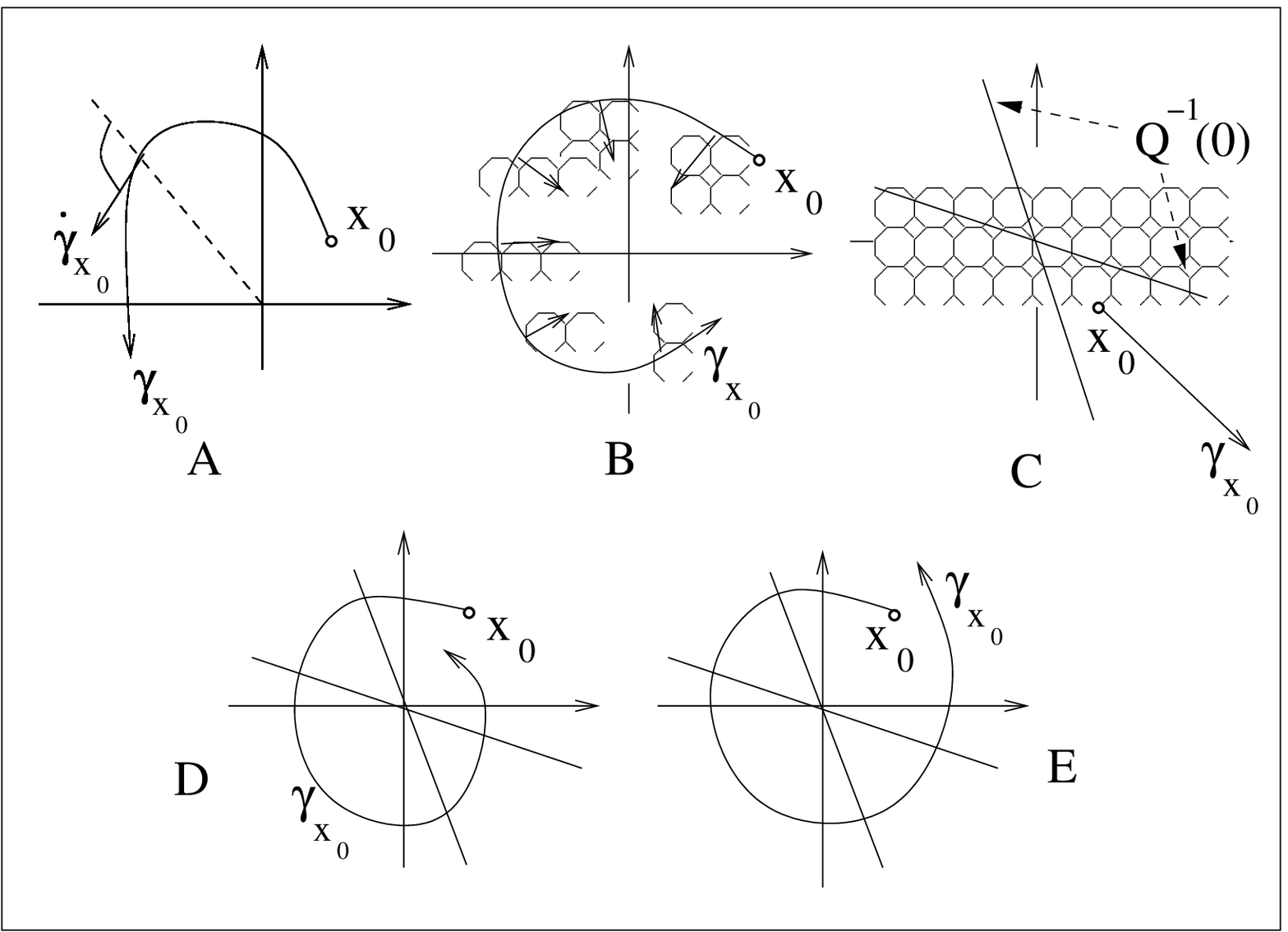}{Proof of the stability conditions}{10}
Theorem \ref{stabb}, stated in Appendix \ref{a-statment}, is the main result of 
\cite{SIAM}, and gives
necessary and sufficient conditions for \pro holding true. We next describe the
main idea of the proof. All  details can be found in \cite{SIAM}.

We build
the ``worst trajectory'' $\g_{x_0}$ i.e. the trajectory (based at $x_0$) having the following property. At each time $t$, $\dot \g_{x_0}(t)$ forms the 
smallest angle (in absolute value) with the
(exiting) radial direction (Figure 1 A). 

Then the system \r{sw-l2b} is \GUES\
if and only if, for each $x_0\in\R^2$, the ``worst trajectory''  $\g_{x_0}$ tends to the origin. The worst 
trajectory is constructed as follows. We study the locus $Q^{-1}(0)$ 
(where $Q(x):=\det(Ax,Bx$)) where the two vector fields $Ax$
and $Bx$ are collinear. The quantity ${\cal D}$, defined in Definition 
\ref{zappa}, is proportional to the discriminant of the quadratic form 
$Q$.
We have several cases:
\bi
\i If $Q^{-1}(0)$ contains only the origin then, in the \cc and \rc
case, one vector field
points
always on the same side of the other and the worst trajectory is a
trajectory of a fixed vector field (either $Ax$ or $Bx$). In that case,
the system is \GUES\ (case {\bf (CC.1)} and {\bf (RC.1)} of Theorem
\ref{stabb}), see Figure 1, case B. 
The situation is similar in case {\bf (RR.1)} (the worst trajectory
tends to the origin). 
\i If $Q^{-1}(0)$  does not contain only the origin then it is the union
of two lines passing through the origin (since $Q$ is a quadratic form). 
If at each point of $Q^{-1}(0)$, the two vector fields have opposite direction,
then there exists a trajectory going to infinity corresponding to a  constant
switching function (see Figure 1, case C). 
This correspond to cases {\bf (CC.2.1)}, {\bf (RC.2.1)}  and {\bf
(RR.2.1)} of Theorem \ref{stabb}. In that situation, there
exists $u\in[0,1]$ such that  the matrix $M(u):=uA+(1-u)B$, 
$u\in[0,1]$ admits an eigenvalue with
positive real part. If at each point of $Q^{-1}(0)$, the two vector fields 
have the same direction, then 
the system is \GUES\ if and only if the worst trajectory turns around
the origin and after one turn the distance from the origin is decreased.
(see Figure 1, cases  D and E).
The quantities $\rocc,\rorc,\rorr$ defined in Theorem \ref{stabb}  
(for the three cases {\bf (CC)}, 
{\bf (RC)}, {\bf (RR)} resp.) represent the distance from 
the origin of the worst trajectory (that at time zero is at distance 
$1$), after one half turn.
This correspond to cases {\bf (CC.2.2)}, {\bf (RC.2.2)} and 
 {\bf (RR.2.2)} of Theorem \ref{stabb}.
\i  Finally {\bf (CC.3)}, {\bf (RC.3)} and {\bf (RR.3)} are the degenerate cases in
which the two straight lines coincide.
\ei
\section{Non Existence of a Uniform Bound on the Minimal Degree of
Polynomial Lyapunov Functions}
\llabel{s-nonexistence}
In this section, we prove Theorem \ref{t-nonexistence}. 
The starting point 
of the argument is to consider a pair of matrices $A$ and $B$ 
having both non real eigenvalues ({\bf (CC)} case) and satisfying: 
\bqn
\D>0,~~\K>1,~~\rocc=1.
\eqnl{azzo}
Such a pair exists. Indeed, Figure \ref{SD1} translates
graphically the contents of Theorem \ref{stabb} for a
fixed $\K>1$, in the region of the $(\roa,\rob)$--plane where
$\roa,\rob>0$.  The open shadowed region
corresponds to values of the parameters $\roa,\rob$ for which the system is GUES. 
We denote by $S^+$ the
open subset of the shadowed region where $\D>0$. The curve $C$ represents
the limit case where $\rocc=1$. 
To each internal point of that 
curve,  it is associated a system verifying 
\r{azzo}, since $\D>0$. 
A system corresponding to such a limit case is not asymptotically stable
but just stable. Moreover, the worst trajectory is a periodic curve,
whose support is of class
$\con^1$ but not of class $\con^2$ (recall that the switchings occur
on $Q^{-1}(0)$, i.e. when the linear vector fields corresponding to $A$ 
and 
$B$ are parallel).

\ppotR{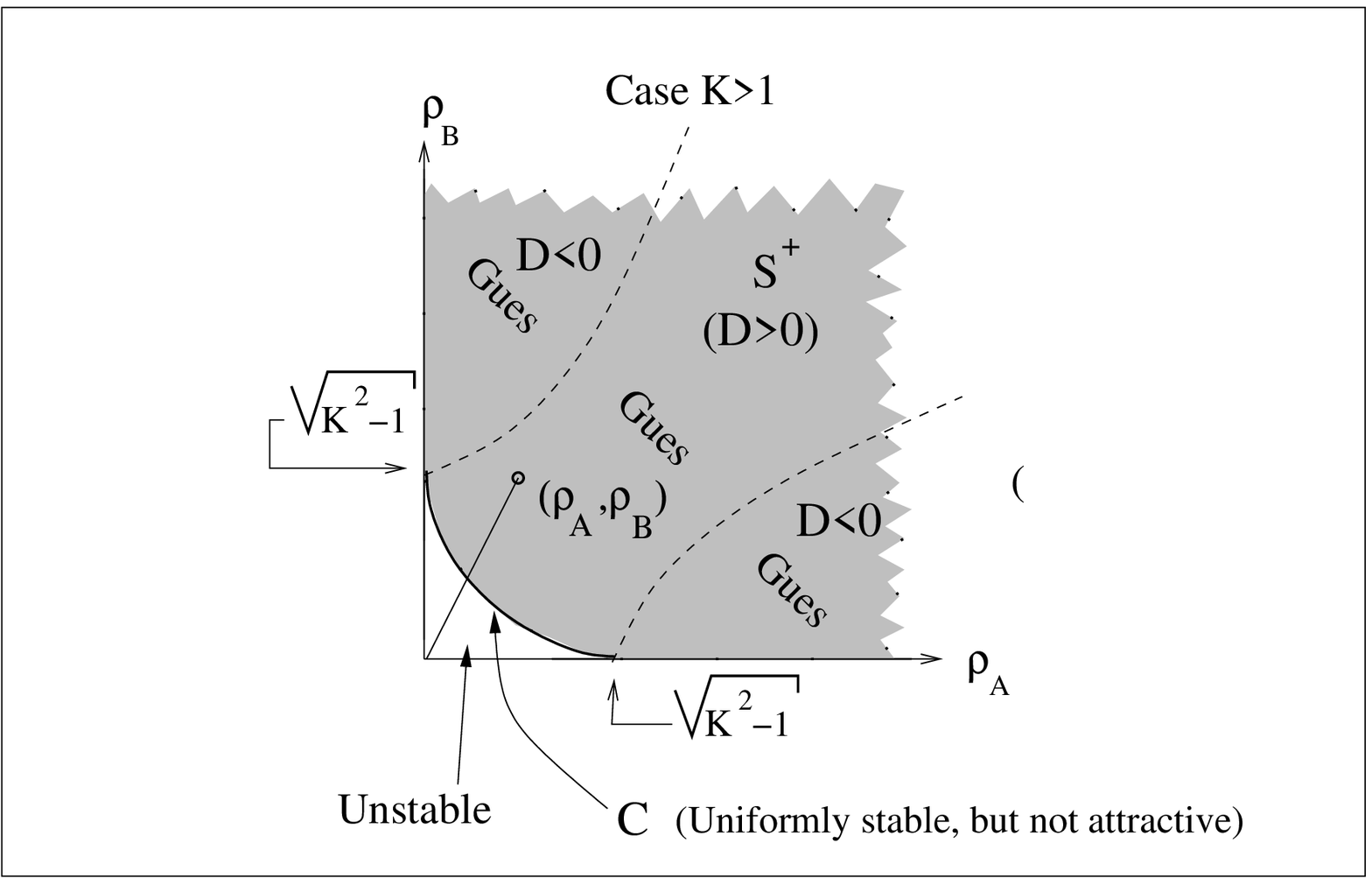}{\GUES\ property in the space of parameters and explicit 
construction of a 5-parameters {\bf SSF} for systems of type 
\r{sw-l2}. $S^+$ is the region in the $(\roa,\rob)$--plane 
in which the system is \GUES\ and $\D>0$.}{11}

Fix a point $(\rho_A,\rho_B)\in C$ corresponding to $(A,B)$.
{}From the picture, it is clear that there exists a sequence of
points $({\roa}_k,{\rob}_k)\in S^+$, for $k\geq 1$,
converging to $(\rho_A,\rho_B)$.  This exactly means that there exists
a sequence of \GUES\ pairs $(A_k,B_k)$, $k\geq 1$, such that $(A_k,B_k)$ tends to 
$(A,B)$ as $k$ goes to $\infty$.

Let $x=(x_1,x_2)$.
For every $k\geq 1$, consider a polynomial CLF 
$V_k=\sum_{1\leq i+j\leq m_k} a_{ij}^{(k)}x^i_1 x^j_2$ of degree at most 
$m_k$.
Arguing by contradiction, we assume that the sequence $(m_k)$
is bounded by a positive integer $m$.
Up to multiplication by a constant, we can choose 
$\sum_{1\leq i+j\leq m_k}|a_{ij}^{(k)}|= 1$. By compactness, there exists a 
subsequence  
of $(V_k)$ (still denoted by $(V_k)$) which converges 
(uniformly on compact subsets of $\mathbb{R}^2$ ) to some non-zero 
polynomial $V$ with degree at most $m$. Note that $V(0)=0$ since the 
$V_k$'s are CLFs.\\

Fix $x_0\in\R^2$, $x_0\not =0$. Let $T>0$ be the period of the worst 
trajectory $\gamma_{x_0}$ corresponding to the pair $(A,B)$, and starting 
at $x_0$ . Note 
that $T$ is independent of $x_0$.
The curve $\gamma_{x_0}:[0,T]\rightarrow \R^2$ can be seen as the 
concatenation of at most five arcs of integral curves of 
$\dot{x}=A\,x$ and $\dot{x}=B\,x\,$ (see Figure \ref{f-nonexistence}) and 
satisfies 
the Cauchy problem:
$$
\left\{ 
\begin{array}{l} 
\dot{x}=C(t)\,x, \\ 
x(0)=x_0,
\end{array} 
\right.
$$
where $C(t)$ is equal to $A$ or $B$ on subintervals of $[0,T]$.\\
For $k\geq 1$, consider the Cauchy problem:
$$
\left\{ 
\begin{array}{l} 
\dot{x}=C_k(t)\,x, \\ 
x(0)=x_0, 
\end{array} 
\right.
$$ 
where $C_k(t)=A_k$ if $C(t)=A$ and $C_k(t)=B_k$ if $C(t)=B$.
Then, $\gamma_k$ is a trajectory of the switched system of the type
\r{sw-l2} associated to $(A_k,B_k)$. 
Since, the right-hand side of the previous equation is Lipschitz continuous in 
$x$ and piecewise continuous in $t$, then the solutions $\gamma_k$ 
converge uniformly to $\gamma_{x_0}$ on $[0,T]$.\\
We next show that $V$ remains constant on $\gamma_{x_0}$. For $k\geq 1$ 
and $t\in [0,T]$,
one has:
$$\|V_k \circ\gamma_k(t)-V\circ\gamma_{x_0}(t)\| \leq \|V_k 
\circ\gamma_k(t)-V\circ\gamma_k(t)\|+\|V\circ\gamma_k(t)-V\circ\gamma_{x_0}(t)\|.$$
By uniform convergence of $V_k$ to $V$ and of $\gamma_k$ to $\gamma_{x_0}$, 
and by continuity of $V$, we deduce that $V_k \circ\gamma_k(t)$ 
converges to $V\circ\gamma_{x_0}(t)$ for every fixed $t$.\\
Since, for every $k\geq 1$, $V_k$ is a CLF for the switched system of the type
\r{sw-l2} associated to $(A_k,B_k)$, then $V_k \circ\gamma_k$ is a decreasing
function and, hence, $V\circ\gamma_{x_0}$ is non-increasing. Moreover
$V\circ\gamma_{x_0}(T)=V\circ\gamma_{x_0}(0)$. Therefore, $V\circ\gamma_{x_0}$ 
must be constant.
It implies that there exists $t_1>0$ such that
either $V(e^{At}x_0)$ or $V(e^{Bt}x_0)$ is constant on $[0,t_1]$. With no loss
of generality, assume the first alternative. Since the map $t\mapsto 
V(e^{At}x_0)$
is real analytic, it follows that $V(e^{At}x_0)$ is constant over the 
whole real line.
By letting $t$ go to $+\infty$, since $e^{At}x_0\to 0$, we deduce that 
$V(x_0)=V(0)=0$. Since $x_0$ is an
arbitrary non zero point of $\R^2$, we get that $V\equiv 0$, which is not possible.

\ppotR{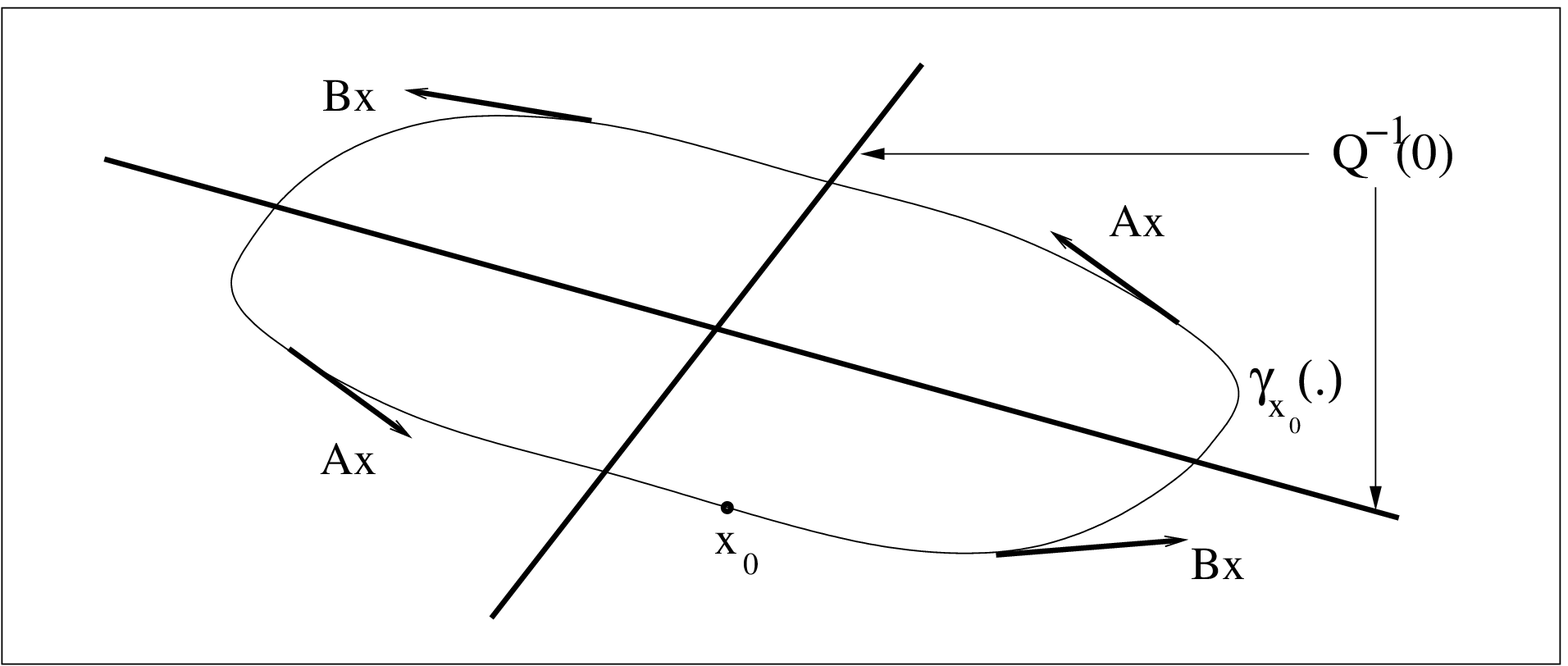}{The ``worst trajectory''}{14}

\brem
The construction of the sequence $(A_i,B_i)$ with unbounded degree 
for polynomial CLF was performed for matrices having both non real 
eigenvalues (that corresponds to the {\bf (CC)} case).
The same construction can be reproduced for the {\bf (RC)} and {\bf (RR)}
cases. 
\erem

\brem 
\llabel{r-ppf-nonexistence}.
For the PPF case (see Remark \ref{r-ppf}), the above argument can be 
easily modified to get  that the minimum total degree of a piecewise 
polynomial CLF
cannot be uniformly bounded over the set of all \GUES\ systems
of the form $\Xi$ (cf. Formula \ref{sw-l2}).
Indeed, let $V_k$ be the sequence of PPFs taking the value  
$V^l_k(x,y)=\sum_{1\leq i+j\leq m}
a_{ijl}^{(k)}x^i_1x^j_2$ in the cone $K^k_{l}$, for $1\leq l\leq
m$. Here, to simplify the notation, we assume without loss of generality,  
that, for each element of the sequence, the number of cones and the 
degree of $V^l_k(x,y)$ is always $m$.

Each cone can be  identified by a couple of angles with the 
$x$-direction. Therefore to each function
$V_k$ we can associate a $m$-uple of angles $(\al^k_1,\ldots,\al^k_m)$ such
that the cone $K^k_{l}$ coincides with the region between the lines 
corresponding to $\al^k_l$ and $\al^k_{l+1}$. In particular, up to 
subsequences, we can assume that the numbers $\al^k_l$ converge to $\al_l$. 
Similarly to the case above, we can normalize the coefficients of the CLFs
$V_k$ by $\sum_{l=1}^m \sum_{1\leq i+j\leq m} |a_{ijl}^{(k)}|=1$ and consider a subsequence
of the coefficients converging to $a_{ijl}$. Then, if we define $V$ as the
PPF such that $V(x,y)=V^l(x,y)= \sum_{1\leq i+j\leq m} a_{ijl} x^i_1x^j_2$ 
on the cone 
$K_l$ defined by the angles $\al_l$ and $\al_{l+1}$, it is easy to verify that 
$V_k(x,y)$ converges uniformly on compact subsets of $\R^2$ to $V(x,y)$.
We can conclude the proof as before showing that $V_l(x_0)=V(0)=0$ for 
arbitrary $x_0$, which leads to a contradiction. 
\erem

\section{Explicit Construction of a finite-{\bf SSF} for Systems 
of Type \r{sw-l2}}
\llabel{s-finite-SSF}
In this section, we provide a 5-parameters finite-{\bf SSF} for the 
class $\csi$ 
of bidimensional systems of type \r{sw-l2}. Recall that 
for what concern the stability issue, $\csi$ can be parameterized  
by the 6-parameters family provided by the pairs  $(A,B)$ (of 
$2\times 2$ matrices) defined up to their norm.

As explained in the introduction, it is enough to construct a CLF for
a pair $(A,B)$ giving rise to a \GUES\ system of $\csi$. We only treat the 
{\bf (CC.2.2)} case since, in all the other cases, the construction is
entirely similar.

In the {\bf (CC)} case, after a three-parameters change of
coordinates, the normal form for the pair $(A,B)$ is given by (cf. 
Appendix \ref{a-normalforms}):
$$
A=\left(\ba{cc}-\rho_A&-1/E\\E&-\rho_A\ea\right),~~~~~~
B=\left(\ba{cc}-\rho_B&-1\\1&-\rho_B\ea\right),~~~~~~
E>0.
$$
Moreover, in the {\bf (CC.2.2)} case, we have $\K>1$, $\D>0$ and
$\rocc<1$, where $\K:= 1/2(E+1/E)$, $\D$ and $\rocc$ being respectively
defined in \r{zappa} and \r{<1cc}. Recall that, for 
fixed 
$\K>1$ (i.e fixed
$E>1$), Figure $2$ describes, in the $(\rho_A,\rho_B)$-plane, the
status of each point with respect to the \GUES\ issue.

We now associate, to every \GUES\ pair $(A,B)$, a pair 
$(\tilde{A},\tilde{B})$ corresponding to a system
of the type \r{sw-l2} uniformly stable but not attractive.
Consider in Figure 2 the line segment joining 
the point $(0,0)$ to $(\rho_A,\rho_B)$ in the $S^+$ region.
That segment intersects the curve $C$ in a point 
$(\tilde{\rho}_A,\tilde{\rho}_B) $. That results from
the Jordan separation theorem and the fact that $C$ connects
the points $(\sqrt{\K^2-1}, 0)$ and  $(0,\sqrt{\K^2-1})$.
Therefore, there exists a $\zeta\in (0,1)$ such that, for
the system given by: 
$$
\tilde{A}=\left(\ba{cc}-\tilde{\rho}_A&-1/E\\E&-\tilde{\rho}_A\ea\right),~~~~~
\tilde{B}=\left(\ba{cc}-\tilde{\rho}_B&-1\\1&-\tilde{\rho}_B\ea\right),~~~~~
\tilde{\rho}_A=\rho_A\zeta,~~~~\tilde{\rho}_B=\rho_B\zeta,
$$
the worst trajectories $\gamma_{x_0}$ are closed curves, i.e. $\rocc=1$.\\
Moreover, one can easily compute: 
$$
\det (\tilde{A}x, Ax)=\ra (1-\zeta)\left({x_1}^2 
E+\frac{{x_2}^2}{E}\right)>0\ ,
$$
$$
\det (\tilde{B}x, Bx)=\rb (1-\zeta)|x|^2>0\ .
$$
Therefore the vector fields $Ax,\,Bx$ point inside the area delimited by a 
fixed worst trajectory (that is closed curve) of the modified switched 
system and so, passing  to angular coordinates,
the function:
\bqn 
V(r,\al)=\frac{r}{\tilde{r}(\al)}\,\label{liav} ,
\eqn
where $\tilde{r}(\al)$ is a parameterization of the fixed worst 
trajectory, is a 
CLF for the system defined by $(A,B)$.\\
Hence, we have provided a 5-parameter {\bf SSF} in the {\bf (CC.2.2)} 
case. The five parameters are: $\K$, the ratio $\rob/\roa$, and   
the three parameters involved in the change of coordinates to get the 
normal forms \r{ARC}, \r{BRC}. 
\brem
Notice that, in the cases {\bf (CC.1)} and {\bf (CC.2.1)} (cf. Section 
\ref{s-bid} and Theorem \ref{stabb} in Appendix \ref{a-statment}),
one can choose as {\bf SSF} the set of quadratic polynomials, which actually 
is parameterized by two parameters.
\erem

\brem
Let us come back to the general system \r{sw-l}, subject to {\bf H0}.
Notice that the question of finding the smallest $m$ such that there exists 
a $m$-parameters finite-{\bf SSF}, for a certain class $C$ of systems,
has no real meaning if one does not require suitable conditions on the 
map $\Psi$ in Definition \ref{d-ssf}. Indeed, it is always possible to build 
a countable {\bf SSF} for the class of systems of 
type \r{sw-l} in $\R^n$ subject to {\bf H0}. 
\erem

\appendix
\section{Stability Conditions for bidimensional Systems}
\llabel{a-A}
\subsection{Statement of the Stability Conditions}
\llabel{a-statment}
\bt
Let $A$ and $B$ be two real matrices such that 
{\bf H1}, {\bf H2},
{\bf H3} and {\bf H4}, given in Section \ref{s-bid},
hold and 
define  $\roa,\rob,\K,\D$ as in Definition
\ref{zappa}.  
We have the following stability conditions:
\bd
\item[Case (CC)] If 
$A$ and $B$
have both complex eigenvalues then:
\bd
\item[Case (CC.1)] if $\D<0$ then \pro is true;                       
\item[Case (CC.2)] if $\D>0$ then:
\bd
\item[Case (CC.2.1)] if $\K<-1$ then \pro is false;
\item[Case (CC.2.2)] if $\K>1$  
then \pro is true if and only if it holds the
following
condition:
\bqn 
\rocc&:=&\exp{\left[
-\roa\,\arctan\left(\frac{-\roa\K+\rob}{\sqrt{\D}}\right) -\right.}\label{<1cc}
\\
&&\hspace{-1cm}
\left.\rob\,\arctan\left(\frac{\roa-\rob\K}{\sqrt{\D}}\right)
-\frac\pi2(\roa+\rob)
              \right]\times\nn\\ 
&&~~~~~\times\sqrt{\frac
{(\roa\rob+\K)+\sqrt{\D}}
{(\roa\rob+\K)-\sqrt{\D}}}<1\nn 
\eqn
\ed
\item[Case (CC.3)] If $\D=0$ then \pro holds true or false whether $\K>1$ or $\K<-1$.
\ed
\item[Case (RC)] If 
$A$ and $B$
have one of them complex and the other real eigenvalues, define
$\chi:=\roa\K-\rob$, where $\roa$ and $\rob$ are chosen in such a way
$\roa\in i\mr,~\rob\in\mr$. Then: \bd
\item[Case (RC.1)] if $\D>0$ then \pro is true;
\item[Case (RC.2)] if $\D<0$ then $\chi\neq0$ and we have:
\bd
\item[Case (RC.2.1)] if $\chi>0$ then $\pro$ is false. Moreover in this
case $\K/i<0$;
\item[Case (RC.2.2)] if $\chi<0$, then: 
\bd
\item[Case (RC2.2.A)] if $\K/i\leq0$  then \pro is true;
\item[Case (RC2.2.B)] if $\K/i>0$  then \pro is true iff
it holds the following condition:
\bqn 
\rorc:={\Big(\frac{m^+}{m^-}\Big)}^{-\frac{1}{2}(\ra/i-1)} e^{-\rb\bar{t}}
\Big(\sqrt{1-\ca^2}\,m^-\,\sin\bar{t}\,-(\cos\bar{t}\,-\frac{\ca}{i}\sin\bar{t}
)\Big)<1\llabel{<1rc}
\eqn

\bqn
\mbox{where:\ \ \ }
&&m^\pm:=\frac{-\chi\pm\sqrt{-\D}}{(-\roa/i-1)\K/i}\nn\\
&&\bar{t}=\displaystyle{\arccos{\frac{-\ra/i+\rb\ca/i}{\sqrt{(1-\ca^2
)(1+\rb^2)}}}}\nn
\eqn
\ed
\ed
\ed
\item[Case (RC.3)] If $\D=0$ then \pro holds true whether
$\chi<0$ or $\chi>0$.
\item[Case (RR)] If $A$ and $B$
have both real eigenvalues then:
\bd
\item[Case (RR.1)] if $\D<0$ then \pro is true. Moreover we have $|\K|>1$;
\item[Case (RR.2)] if $\D>0$ then $\K\neq -\roa\rob$ (notice that
$-\roa\rob>1$) and :
\bd
\item[Case (RR.2.1)] if $\K>-\roa\rob$ then $\pro$ is false 
\item[Case (RR.2.2)] if $\K<-\roa\rob$ then:
\bd
\item[Case (RR.2.2.A)] if $\K>-1$ then $\pro$ is true;
\item[Case (RR.2.2.B)] if $\K<-1$ then $\pro$ is true iff the following
condition holds:
\bqn
&&\rorr:=-f^{sym}(\roa,\rob,\K)f^{asym}(\roa,\rob,\K)\times\label{<1rr}\\
&&~~~~~~~~~~~~~~~f^{asym}(\rob,\roa,\K)<1,
\nn
\eqn
where:
\bqn
f^{sym}(\roa,\rob,\K)&:=&
\frac{1+\roa/i+\rob/i+\K-\sqrt{\D}}{1+\roa/i+\rob/i+\K+\sqrt{\D}};\nn\\
f^{aym}(\roa,\rob,\K)&:=&\left(
\frac{\rob/i-\K\roa/i-\sqrt{\D}}{\rob/i-\K\roa/i+\sqrt{\D}}\right)^
{\frac12(\roa/i-1)}.
\eqnn
\ed
\ed
\item[Case (RR.3)] If $\D=0$ then \pro holds true or false 
whether $\K<-\roa\rob$ or $\K>-\roa\rob$.
\ed
\ed
Finally, if \pro does not hold true, then in case {\bf CC.2.2} with $\rocc=1$, case
{\bf
(RC.2.2.B)},  with $\rorc=1$, case
{\bf
(RR.2.2.B)},   with $\rorr=1$,
case {\bf (CC.3)} with $\K<-1$, case 
{\bf (RC.3)} with $\chi>0$ and 
case {\bf (RR.3)} with $\K>-\roa\rob$,
the origin is just stable.
In
the other cases, the system is unstable.
\llabel{stabb}
\et
\brem
\llabel{r-correct}
Formula \r{<1rr} is a corrected version of Formula (6), p.93, of \cite{SIAM} and 
it is proved in Appendix \ref{a-proof}.
\erem
\subsection{Normal Forms of $2\times 2$ Matrices}
\llabel{a-normalforms}
\bp\llabel{p-normal}
Let $A$, $B$ be two $2\times2$ real matrices satisfying  conditions {\bf H1}, 
{\bf H2}, {\bf H3} and {\bf H4} given in Section \ref{s-bid}. 
In the case in which one of the two matrices
has real and the other non-real eigenvalues (i.e. the \rc case),
assume that $A$ is the one having real eigenvalues.
Then there 
exists a 3-parameter change of coordinates and two 
constant $\al_A,\al_B>0$ such that the matrices 
$A/\al_A$ and $B/\al_B$ (still denoted below by $A$ and $B$) 
are in the following normal forms:
\bd
\i[Case in which $A$ and $B$ have both non-real eigenvalues (\cc case):]
\bqn
A=\left(\ba{cc}-\roa&-1/E\\E&-\roa\ea\right),~~~~~~~~
B=\left(\ba{cc}-\rob&-1\\1&-\rob\ea\right),
\eqn
where $\roa,\rob>0$,  $|E|>1$. In this case,
$\K=\frac{1}{2}(E+\frac{1}{E})$. 
Moreover, the 
eigenvalues of $A$ and $B$ are respectively 
$-\roa\pm i$ and $-\rob\pm i$.
\i[Case in which $A$ has real and $B$ non-real eigenvalues (\rc case):]
\bqn
&&A=\left(\ba{cc}-\roa/i+1&0\\0&-\roa/i-1\ea\right),\llabel{ARC}\\
&&B=
\left(\begin{array}{c c}
-\rb-\ca/i & -\sqrt{1-\ca^2} \\
\sqrt{1-\ca^2} & -\rb+\ca/i
\end{array}\right),\llabel{BRC}
\eqn
where $\rob>0$, $\roa/i>1$, $\K\in i\R$.  In this case, the 
eigenvalues of $A$ and $B$ are respectively $-\roa/i\pm1$ and $-\rob\pm 
i$.
\i[Case in which $A$ and $B$ have both real eigenvalues (\rr case):]
\bqn
A&=&\left(\ba{cc}   -\roa/i+1   &0\\0&  -\roa/i-1\ea\right),\label{ARR}\\
B&=&
\left(\ba{cc}  \K-\rob/i&   1-\K\\ \
1+\K& -\K-\rob/i\ea\right),\label{BRR}\\
\eqn
where $\roa/i,\rob/i>1$ and $\K\in\mr\setminus\{\pm1\}$. In this case, the
eigenvalues of $A$ and $B$ are respectively $-\roa/i \pm1 $  and
$-\rob/i\pm1$.
\ed
\ep
\subsection{Proof of Formula \r{<1rc}}\llabel{a-proof}
In this paragraph, we prove Formula \r{<1rc}, i.e. in the 
\textbf{(RC.2.2.B)} case, we determine 
an inequality defining the set of parameters $\ra,\ \rb,\ \K $ such 
that the property
\pro, stated in Section \ref{ss-bid}, holds.

Thanks to Proposition \ref{p-normal} (see also \cite{SIAM}, Appendix B, 
p.110), we can find a coordinate 
transformation such that (up 
to a rescaling of the matrices) $A$ and $B$ are given by 
equations \r{ARC}, \r{BRC}.
In the case \textbf{(RC.2.2.B)}, we have 
$\Dd:=\K^2+2\roa\rob\K-
(1+\roa^2+\rob^2)<0$, $\chi:=\roa\K-\rob<0$, $\ca/i>0$.
Moreover, the set $Q^{-1}(0)$ is the union of two
lines passing from the origin and, at each point of $Q^{-1}(0)$, the two 
vector fields point in the same direction.
One easily checks that the slope of the two 
lines defining $Q^{-1}(0)$
is:
$$
m^{\pm}=\frac{-\chi\pm\sqrt{-\Dd}}{(-\ra/i-1)
\sqrt{1-\ca^2}}.
$$ 
Notice that, in our case we have $m^{\pm}<0$ and $m^+<m^-$.

In this case, the worst trajectories are concatenations of arcs 
of integral curves of the vector fields $Ax$, $Bx$ and rotate counterclockwise 
around the origin. More precisely, they are integral curves of $Ax$ 
from the line $x_2=m^+ x_1$ to the line $x_2=m^- x_1$, and integral 
curves  of $Bx$ otherwise.

Therefore, starting from the point $\vett{1}{m^+}$ (with the field $Ax$), we 
follow the worst trajectory until it touches again the line $x_2=m^+ x_1$. 
Property \pro is then satisfied if and only if $\rorc<1$, where 
$\rorc$ is the absolute value of the first coordinate of the final point. 

One can easily compute that the first switching time is 
$t_1=\displaystyle{\frac{1}{2}
\log{\frac{m^+}{m^-}}}\,$, which is positive since $\displaystyle{\frac{m^+}{m^-}>1}
$. Moreover, the integral curve of $Bx$ starting from the point 
$\vett{1}{m^-}$ is:
$$
e^{-\rb t}\vett{-\sqrt{1-\ca^2}\,m^-\,\sin t\,+  
(\cos t\,-\frac{\ca}{i}\sin t)}{\sqrt{1-\ca^2}\,\sin t\,+
\,m^-(\cos t\,+\frac{\ca}{i}\sin t)},
$$
and, setting the ratio between the second 
coordinate and the first one equal to $m^+$, one obtains that the second 
switching 
time is $t_2=\displaystyle{\arccos{\frac{-\ra/i+\rb\ca/i}{\sqrt{(1-\ca^2
)(1+\rb^2)}}}}$.
Notice that $t_2$ is well defined if and only if $\Dd<0$ (condition which is 
satisfied in our case). Moreover, $t_2$ is positive and less than 
$\pi$. Finally, the inequality we was seeking 
for is:
$$\rorc={\Big(\frac{m^+}{m^-}\Big)}^{-\frac{1}{2}(\ra/i-1)} e^{-\rb t_2}
\Big(\sqrt{1-\ca^2}\,m^-\,\sin t_2\,-(\cos t_2\,-\frac{\ca}{i}\sin t_2
)\Big)<1\,.$$

\section{The Stability Properties of \r{sw-l} Only Depend on the 
Convex Hull of $\aaa$}
\llabel{s-b}
We provide here the proof of Proposition~\ref{p-convex}.
First, let us show the following:\\\\
{\bf Claim} Consider the switched system \r{sw-l},  
under {\bf H0}, and let $\aaa'$ be a 
measurable subset of $\aaa$ such that the convex
hull of $\aaa'$ contains $\aaa$. 
Then the following two conditions are equivalent:
\bd 
\i[i)] the system is \GUES\ (resp. uniformly stable), 
with $A(.)$ measurable, taking values in $\aaa$,

\i[ii)] the system is \GUES\ (resp. uniformly stable),
with $A(.)$ measurable, taking values in $\aaa'$.
\ed
\noindent
{\bf Proof of the Claim.}
Let $\Aa$ (resp. $\Aa'$) be the set of measurable functions 
$A(.):[0,\infty[\to \aaa$ (resp. $A(.):[0,\infty[\to \aaa'$).
Since $\aaa'$ is contained in $\aaa$ then the implication {\bf i)} 
$\Rightarrow$ {\bf ii)} is obvious.

Let us prove the other implication (that is strictly related to the classical 
approximability theorems in control theory). We start considering 
uniform stability.
By contradiction, assume that we can find $\epsilon>0$ satisfying the 
following.  There exists a sequence of points $(x_l)$ tending to zero
and a sequence of controls $A_l(.)\in \Aa$ such that the 
corresponding trajectory $\gamma_l$ starting at $x_l$ exits 
the interior of the ball of radius $\epsilon$ for some time $t_l$.
Using classical 
approximability results (see for instance \cite{agra-book}),  
the trajectory $\gamma_l$ can be approximated in the $L^\infty$-norm 
on $[0,t_l]$ by a
trajectory $\gamma'_l$ corresponding to a switching function
${A_l'}(.)\in \Aa'$ and starting at $x_l$. Hence $\gamma'_l$ exits the
interior of the ball of 
radius $\epsilon/2$ at time $t_l$. We reached a contradiction.

Now we want to prove that \GUES\ holds in the case  
$A(.)\in \Aa'$ implies \GUES\ holds in the case $A(.)\in \Aa$.
Since $\aaa$ is compact, we 
know (see Definition \ref{d-UAS} and below) 
that attractivity and \GUES\ are equivalent for the
 corresponding switched system.
Therefore, proceeding by contradiction, we can assume that there is a 
trajectory
$\gamma(.)$ of the switched system corresponding to $A(.)\in \Aa$ 
not converging to zero. That means that there exist $\epsilon>0$ and a 
sequence $t_n$ of times tending to infinity such that 
$|\g(t_n)|>\epsilon$.
As before, we can approximate $\gamma(.)$ on the interval $[0,t_n]$ with a 
trajectory $\gamma_n(.)$ corresponding to controls taking values in 
$\aaa'$, in such a way that $|\gamma_n(t_n)|>\epsilon/2$ . 
But this is impossible since we have assumed \GUES\ for the switched 
system with $A(.)\in \Aa'$. 

Notice that one can provide an alternative argument for the \GUES\ part
of Proposition~\ref{p-convex}, by using CLFs.
\hspace{\stretch{1}}$\Box$

Then one immediately extend to semicones observing that the stability 
properties of the system \r{sw-l}, subject to the compactness hypothesis 
{\bf H0}, depend only on the shape of the trajectories and not on the way 
in which they are parameterized.

\end{document}